\documentclass[a4,11pt,epsf, amssymb]{article}
\usepackage{color}
\usepackage{graphicx}
\usepackage{hyperref}
\usepackage{amsmath,amsfonts,latexsym,amssymb}
\begin{document}
%\voffset = +18mm
%\renewcommand{\theequation}{\arabic{section}.\arabic{equation}}
\renewcommand{\theequation}{\arabic{equation}}
\title{\bf A Note on The Backfitting Estimation of Additive Models}
\author{\\ \normalsize
Yingcun  Xia
\\
\
Guizhou College of Finance and Economics, China
\\
National University of Singapore, Singapore\\
Email: \textit{staxyc@nus.edu.sg}  }
\date{}
\def\beginn{\begin{eqnarray*}}
\def\endn{\end{eqnarray*}}
\def\beginy{\begin{eqnarray}}
\def\endy{\end{eqnarray}}
\def\n{\nonumber}
\newtheorem{Theorem}{Theorem}[section]
\newtheorem{Example}[Theorem]{Example}
\newtheorem{Lemma}[Theorem]{Lemma}
\newtheorem{Note}[Theorem]{Note}
\newtheorem{Proposition}[Theorem]{Proposition}
\newtheorem{Corollary}[Theorem]{Corollary}
\newtheorem{Remark}[Theorem]{Remark}
\def\ditem{\vspace{-0.1cm} \item}

\def\Cov{\mbox{\rm Cov}}
\def\Var{\mbox{\rm Var}}
\def\YX{_{Y|X}}
\def\argmin{\mbox{\rm arg}\min}
\def\b{{\mathbf b}}
\def\pperp{\perp\hspace{-.25cm}\perp}
\def\t{{\hspace{-0.05cm}\top}}
\def\Bt{B^\top}
\def\E{{\cal E}}
\def\D{{\cal D}}
\def\N{{\cal N}}
\def\O{{\cal O}}
\def\V{{{\cal V}}}
\def\R{{\Bbb R}}
\def\btd{\bigtriangledown}
\def\btdt{\bigtriangledown^\t\hspace{-0.1cm}}
\def\dfor{{\qquad \mbox{\rm for} \quad}}
\def\C{{\cal C}}
\def\rt{\raisebox{-1.5ex}[0pt]}
\def\bbar{{b\hspace{-0.1cm}\bar{} \ }}
\def\diag{\mbox{\rm diag}}
\def\tg{\tilde g}
\def\eqc{&\hspace{-0.35cm}=&\hspace{-0.35cm}}
\def\deqc{&\hspace{-0.35cm}\stackrel{def}=&\hspace{-0.35cm}}
\def\S{\mathbf{S}}

\maketitle

\renewcommand{\theTheorem}{\arabic{Theorem}}

\rm

\begin{abstract}

\baselineskip2.2em

{\small The additive model is one of the most popular semiparametric
models. The backfitting estimation (Buja, Hastie and Tibshirani,
1989, \textit{Ann. Statist.}) for the model is intuitively easy to understand and
theoretically most efficient (Opsomer and Ruppert, 1997,
\textit{Ann. Statist.}); its implementation is equivalent to solving
simple linear equations. However, convergence of the
algorithm is very difficult to investigate and is still unsolved. For
bivariate additive models, Opsomer and Ruppert (1997, \textit{Ann. Statist.}) proved the convergence under a very strong condition
and conjectured that a much weaker condition is sufficient.  In this short note, we show that a weak condition can guarantee the convergence of the backfitting estimation algorithm when the Nadaraya-Watson kernel smoothing is used.}

\vspace{.5cm}


\end{abstract}

%\newpage

%\pagebreak
\vspace{.7cm}

\baselineskip2.0em

\setcounter{equation}{0}

\section{Introduction}

The additive model has been proved to be a very useful
semiparametric model and is popularly used in practice. An intuitive
implementation of the estimation is the backfitting approach (Buja,
Hastie and Tibshirani, 1989, called BHT hereafter). It is noticed that the implementation can be done easily by solving linear normal equations (pp. 476, BHT) if \textit{the backfitting algorithm converges.} However, to justify the convergence of the
algorithm is not easy.  BHT provided sufficient conditions that
guarantee the convergence of the backfitting algorithm or,
equivalently, the existence of the estimators. These conditions are
only generally satisfied by regression splines and other methods, but not by kernel smoothing. Some other approaches (e.g.
Tj{\o}stheim and Auestad, 1994; Linton and Nielsen, 1995; Mammen,
Linton and Nielsen,  1999; Wang and Yang, 2007) have been proposed
to avoid hard problems about the convergence of algorithm and the
asymptotics of estimators. However, the original backfitting of BHT
is still one of the most intuitive approach.

Opsomer and Ruppert (1997, called OR hereafter) investigated the algorithm's convergence
for the local polynomial kernel smoothing when the predictors are bivariate. Suppose $ Y $ is the
response and $ (U, V) $ is the bivariate predictors satisfying
the additive model
\begin{equation}
 Y = \alpha + m_1(U) + m_2(V) + \varepsilon, \label{model1}
\end{equation}
where $ E(\varepsilon|U, V) = 0 $ almost surely. Constraints $ E \{m_1(U)\} = E\{ m_2(V)\} = 0$ are usually imposed for  model identification; see for example OR. It is known (see, e.g. BHT) that the terms in the model are the solution to
minimizing
\begin{equation}
  \min_{\stackrel{m_1 \in L_2, m_2 \in L_2,}{_{\alpha \in R}}} E\{ Y - \alpha - m_1(U) -
  m_2(V) \}^2,
\end{equation}
where $ L_2 $ is the measurable functional space with finite second moments. Let $
f(u, v) $, $f_1(u) $ and $ f_2(v) $ be the joint density
function and marginal density functions of $ (U, V) $, $ U $
and $ V $ respectively. OR required that
  $$
  \sup_{u, v } \Big|\frac{f(u, v)}{ f_1(u) f_2(v)} - 1 \Big| < 1
  $$
to prove the convergence of the backfitting algorithm. This
requirement is very stringent and even excludes a big part  of the normal
distributions. However, OR conjectured that the algorithm convergence can be
guaranteed under very week conditions. Next, we shall prove that their conjecture is correct when the Nadaraya-Watson kernel is used.

\section{Main results}
Suppose $ \{(Y_i, U_{i}, V_{i}): i = 1, ..., n\} $ is a random sample
from  model (\ref{model1}).  Following BHT, let $ \mathbf{m}_1 =
(m_1(U_{i}), ..., m_1(U_{n}))^\top$,  $ \mathbf{m}_2 =
(m_2(V_{i}), ..., m_2(V_{n}))^\top$ and $ \mathbf{Y} = (Y_1,
..., Y_n)^\top $. The estimators of functions $ m_1 $ and $ m_2 $ are determined by the estimation of funcation values at the observed points, i.e. $ \mathbf{m}_1 $ and $ \mathbf{m}_2$.  Let $ K(.) \ge 0 $ be kernel function and $ K_{h}(.) = K(./h)/h$ for any $ h>0$.

For the estimation of function values at $ U_i$ and $ V_i $, we use (varying) bandwidth $ h_i > 0 $ and $
\hbar_i >0 $ respectively and kernel weights
 $ \ell_i = [K_{h_i}(U_i-U_1), ...,
K_{h_i}(U_i- U_n)]^\top/\sum_{k=1}^n K_{h_i}(U_i-U_k)
$ and  $ \omega_i =
[K_{\hbar_i}(V_i-V_1), ...,
K_{\hbar_i}(V_i- V_n)]^\top/\sum_{k=1}^n K_{\hbar_i}(V_i-V_k)
 $. Let
$$
 \S_1 = \begin{pmatrix} \ell_1^\top \\
     \vdots\\
     \ell_n^\top\end{pmatrix}, \qquad
 \S_2 = \begin{pmatrix} \omega_1^\top \\
     \vdots\\
     \omega_n^\top\end{pmatrix}     .
$$
Corresponding to constraints $ E\{m_1(U)\} = E\{m_2(V)\} = 0$, we introduce $(\mathbf{I}_n - \mathbf{1}_n \mathbf{1}_n^\top/n) $, where $
\mathbf{I}_n $ is the $ n\times n $ identity matrix and $ \mathbf{1}_n
$ is a vector of $ n\times 1 $ with all entries 1. Let $ \S_1^*
=(\mathbf{I}_n - \mathbf{1}_n \mathbf{1}_n^\top/n) \S_1 $ and $ \S_2^*
=(\mathbf{I}_n - \mathbf{1}_n \mathbf{1}_n^\top/n) \S_2 $. Using kernel smoothing, the
backfitting estimation procedure is iteratively
$$
  \hat{\mathbf{m}}_1^{new} := \S_1^*\{ \mathbf{Y} -\hat{\mathbf{ m}}_2^{old}\}, \qquad \hat{\mathbf{m}}_2^{new} := \S_2^*\{ \mathbf{Y} -
  \hat{\mathbf{m}}_1^{old}\}.
$$
As BHT pointed out, the final estimators $\hat{\mathbf{m}}_1 $ and $ \hat{\mathbf{m}}_2 $ of the algorithm are equivalent to the solution of
\beginn
    \begin{pmatrix} \mathbf{I}_n & \S_1^* \\ \S_2^* & \mathbf{I}_n \end{pmatrix}
    \begin{pmatrix}\hat{\mathbf{m}}_1 \\ \hat{\mathbf{m}}_2 \end{pmatrix}
     = \begin{pmatrix} \S_1^* \\ \S_2^* \end{pmatrix} \mathbf{Y}.
\endn
The solution exists if the inverse  of $ (\mathbf{I}_n
- \S_2^* \S_1^*) $ or $ (\mathbf{I}_n - \S_1^* \S_2^*) $ exits. If the iteration converges, then
estimators of $ \alpha $, $ \hat{\mathbf{m}}_1 $ and $\hat{\mathbf{m}}_2$ are respectively $ \hat \alpha = \bar Y $,
\begin{equation*}
% \hat{\mathbf{m}}_1 = \{\mathbf{I}_n- (\mathbf{I}_n - \S_1^{*} \S_2^{*})^{-1}(\mathbf{I}_n-\S_1^{*})\}  \mathbf{Y}
 \hat{\mathbf{m}}_1 = \S_1^*(\mathbf{I}_n - \S_2^{*} \S_1^{*})^{-1}(\mathbf{I}_n-\S_2^{*})  \mathbf{Y}
\end{equation*}
and
\begin{equation*}
% \hat{\mathbf{m}}_2 = \{\mathbf{I}_n- (\mathbf{I}_n - \S_2^{*} \S_1^{*})^{-1}(\mathbf{I}_n-\S_2^{*})\}  \mathbf{Y}
 \hat{\mathbf{m}}_2 = (\mathbf{I}_n - \S_2^{*} \S_1^{*})^{-1}\S_2^*(\mathbf{I}_n-\S_1^{*})  \mathbf{Y}
\end{equation*}
(the solutions can be rewritten in different forms). As we can see, the backfitting estimation is very easy to implement and is equivalent to a one-step calculation, if it converges. Thus, convergence of the algorithm is  essential for the estimation of the additive model.

\begin{Theorem} Denote the order statistics of $ \{U_{1}, ..., U_n\}$  and $
\{V_1, ..., V_n\} $ by $ \{U_{[1]}, ..., U_{[n]}\}$ and $\{V_{[1]}, ..., V_{[n]}\} $ respectively, and their corresponding bandwidths by $ \{h_{[1]}, ..., h_{[n]}\}$ and $\{\hbar_{[1]}, ..., $ $ \hbar_{[n]}\} $ respectively.  If kernel function $
K(.) $ and the bandwidths satisfy $ K(0)>0 $,
\beginy
  \begin{array}{l}   K_{h_{[i]}} (U_{[i]} - U_{[i-1]}) >0 , \qquad K_{h_{[i]}} (U_{[i]} - U_{[i+1]})
    >0,\\
    \vspace{-0.2cm}
    \\
    K_{\hbar_{[i]}} (V_{[i]} - V_{[i-1]}) >0 , \qquad \ K_{\hbar_{[i]}} (V_{[i]} - V_{[i+1]})
        >0,
    \end{array} \label{condition1}
\endy
for $ 1< i < n $, and
\beginy
  \begin{array}{l}  K_{h_{[1]}} (U_{[1]} - U_{[2]}) >0 , \qquad K_{h_{[n]}} (U_{[n]} - U_{[n-1]})
    >0,\\
    \vspace{-0.2cm}
    \\
    K_{\hbar_{[1]}} (V_{[1]} - V_{[2]}) >0 , \qquad \ K_{\hbar_{[n]}} (V_{[n]} - V_{[n-1]})
    >0,    \end{array}  \label{condition2}
\endy
 then the backfitting algorithm converges.
\end{Theorem}

\setcounter{Theorem}{0}

\begin{Remark} Suppose $ K(.) $ is a symmetric kernel function with $ K(v)>0 $ for all $ |v| < 1$
and that global (constant)  bandwidthes $ h $ and $ \hbar $ are used. If $ h $
and $ \hbar $ are bigger than the largest difference between any two
nearest points respectively, i.e.
\begin{equation}
h>
 \max\{ U_{[i+1]} - U_{[i]}, i = 1, ..., n-1 \} \ \mbox{ and } \ \hbar >  \max\{ V_{[i+1]} - V_{[i]}, i = 1, ..., n-1\},   \label{cc}
\end{equation}
then (\ref{condition1}) and (\ref{condition2}) hold. By Theorem 1 the convergence of the algorithm is guaranteed.
\end{Remark}

\setcounter{Theorem}{0}
\begin{Corollary}
Suppose $ U $ and $ V $ are distributed on two compact intervals respectively with density functions bounded away from 0. If global (constant) bandwidths   $ h  $ and $ \hbar  $ are used with $  h, \hbar \to 0 $ and $ nh/\log(n),  n\hbar/\log(n) \to \infty$, then the algorithm converges in probability as $ n $ is large enough.
\end{Corollary}

\begin{Remark} It is remarkable that the range of bandwidths for the algorithm to converge is quite wide, and that bandwidths $ h \propto n^{-\delta} $ and $ \hbar \propto n^{-\delta} $ with $0< \delta < 1 $ satisfy the requirement in Corollary 1. Thus, the algorithm converges. These bandwidths include the optimal bandwidths where $ \delta = 1/5 $ (see, e.g. OR).
\end{Remark}

This short note only considers the bivariate case with  Nadaraya-Watson kernel smoothing. We conjecture that the backfitting estimation still converges under weak conditions for  general additive models and other kernel estimation methods including the local polynomial smoothing. After the convergence is justified, asymptotics of the estimators can be obtained following exactly the same arguments of Opsomer  and Ruppert (1997). The details are omitted.

\section{Proofs}

The proof of Theorem 1 is based on the properties of the regular Markov chain and the Perron-Frobenius theorem (see, e.g. Minc, 1988). The proof of Corollary 1 is based on the properties of order statistics (see, e.g. David and Nagaraja, 2003).

\textbf{Proof of Theorem 1}. We first prove that the absolute  eigenvalues of $ \S_1 $ are all smaller than 1 with only one exception that equals 1. It is easy to see that $ \S_1 $ is a probability transition matrix  of the Markov chain. By conditions (\ref{condition1}) and (\ref{condition2}), $ \S_1 $ is irreducible and aperiodic. Therefore it is a regular transition probability matrix. There is an integer $ k $ such that all entries in $ \S_1^k $ are strictly positive (see, e.g. Romanovsky, 1970, Theorem 14.I).   By the Perron-Frobenius theorem, there is one (and only one) eigenvalue $ \lambda_1 $ of multiplicity 1 such that all entries in its corresponding eigenvector are positive.  It is easy to see that this eigenvalue is $ \lambda_1 = 1$ and its eigenvector is $ \theta = \mathbf{1}_n/\sqrt{n} $, because the sum of any row in $ \S_1 $ is 1. Let $ \lambda_2, ..., \lambda_n $ be the other $n-1 $ eigenvalues of $ \S_1$ (repeated eigenvalues are counted repeatedly). The Perron-Frobenius theorem also indicates that $ 1 = \lambda_1 > \max\{ |\lambda_2|, ..., |\lambda_n| \} $.

Next, we show that the absolute eigenvalues of  $ \S_1^* = (\mathbf{I}_n-\theta\theta^\top) \S_1 $ are all strictly smaller than 1. Suppose that the eigenvalues $ \lambda_2, ..., \lambda_n $ of $ \textbf{S}_1 $ are distinct and their corresponding eigenvectors are $ \beta_2, ..., \beta_n $ respectively (The general argument is similar, but needs more complicated notation).  It is easy to check that $ \theta $ and $ (\mathbf{I}_n-\theta\theta^\top)\beta_k, k = 2, ..., n $ are the eigenvectors of $\S_1^* $ with corresponding eigenvalues being $ 0 $ and $ \lambda_2, ..., \lambda_n $ respectively, because
\beginn
 (\mathbf{I}_n-\theta\theta^\top) \S_1 \theta = (\mathbf{I}_n-\theta\theta^\top) \lambda_1 \theta = 0
\endn
and
\beginn
 (\mathbf{I}_n-\theta\theta^\top) \S_1 (\mathbf{I}_n-\theta\theta^\top) \beta_k &=& (\mathbf{I}_n-\theta\theta^\top) \{ \S_1 \beta_k- \S_1\theta\theta^\top\beta_k\} \\
    &=& (\mathbf{I}_n-\theta\theta^\top) \{ \lambda_k \beta_k- \lambda_1 \theta\theta^\top\beta_k\} \\
    &=&  \lambda_k (\mathbf{I}_n-\theta\theta^\top) \beta_k - \lambda_1 (\mathbf{I}_n-\theta\theta^\top)\theta\theta^\top\beta_k \\
    &=&  \lambda_k (\mathbf{I}_n-\theta\theta^\top) \beta_k, \quad
    \mbox{ for } k = 2, ..., n.
\endn
Since the absolute values of $ 0, \lambda_2, ..., \lambda_n $ are all smaller than 1,
we proved that the absolute  eigenvalues of $
\S_1^* $ are smaller than 1. Applying the same argument to $ \S_2^* $, we have the absolute values of all
eigenvalues of $ \S_2^* $ are smaller than 1.

Since the largest absolute eigenvalues of both $ \S_1^* $ and $ \S_2^* $ are smaller than 1, the absolute values of all eigenvalues of $ \S_2^* \S_1^* $ and $ \S_1^* \S_2^* $  are also smaller than 1. It follows that the inverses of $
(\mathbf{I}_n - \S_2^* \S_1^*) $  and $ (\mathbf{I}_n - \S_1^* \S_2^*) $ exist, and thus the algorithm converges. $\hspace{\fill}\Box$

\textbf{Proof of Corollary 1}.  It is easy to check
\beginy
&&  P(h>
\n \max\{ U_{[i+1]} - U_{[i]}, i = 1, ..., n-1 \}, \  \hbar >  \max\{ V_{[i+1]} - V_{[i]}, i = 1, ..., n-1\}) \\
\n &&  \ge 1 - P( h \le \max\{ U_{[i+1]} - U_{[i]}, i = 1, ..., n-1 \}) \\
&& \hspace{6cm} - P(\hbar \le  \max\{ V_{[i+1]} - V_{[i]}, i = 1, ..., n-1\}).  \label{kdf}
\endy
Consider the second term above. We have
\begin{equation}
 P( h \le \max\{ U_{[i+1]} - U_{[i]}, i = 1, ..., n-1 \}) \ \le \ \sum_{i=1}^{n-1} P( h \le U_{[i+1]} - U_{[i]}).  \label{kkk}
\end{equation}
Let $ F$ be the cumulative probability function of $ U $. Then $ U'= F(U) $ is uniformly distributed on [0, 1]. Let $ U'_{[i]} = F(U_{[i]}) $. By the joint distribution of $(U'_{[i]}, U'_{[i+1]}) $ (see, e.g. David and Nagaraja, 2003) and simple calculation, we have
for any $ c> 0$
\beginn
  P( c \le U'_{[i+1]} - U'_{[i]}) &=& \int_{ \tilde u > u + c} \frac{n!}{(i-1)! (n-i-1)!} u^{i-1} (1-\tilde u)^{n-i-1}  du d\tilde u  \\
    &=& \left\{\begin{array}{ll} (1-c)^n, & \mbox{if }  0\le c \le 1, \\ 0, & \mbox{if } c>1.
    \end{array}\right.
\endn
Let $ c_0 = \inf\{ f^{-1}(u),  0 \le u \le 1\} $, which is positive by the assumption. Note that $ U_{[i]} = G(U'_{[i]}) $, where $ G $ is the inverse function of $ F $. By the property of inverse function, we have $ U_{[i+1]} - U_{[i]} \le c_0  (U'_{[i+1]} - U'_{[i]}) $. Thus
$$
 P( c_0 c  \le U_{[i+1]} - U_{[i]}) \ < \ P(  c  \le U'_{[i+1]} - U'_{[i]}) \ = \   \left\{\begin{array}{ll} (1-c)^{n-1}, & \mbox{if }  0\le c \le 1, \\ 0, & \mbox{if } c>1.
    \end{array}\right.
$$
When $ n $ is large, we can assume $ h < 1 $. It follows that
\beginy
\n \sum_{i=1}^{n-1} P( h \le U_{[i+1]} - U_{[i]}) &\le& n (1-h)^{n-1} \ = \ n \exp\{ (n-1)\log(1-h)\} \\
\n &\le& n\exp\{( n-1) ( -h + h^2/2)\} \ \le \ n\exp\{ -(n-1)h/2\}\\
& \to& 0  \label{jj}
\endy
as $ n \to \infty $. Condition $ n h/\log(n) \to \infty $ is used in the last step of (\ref{jj}). By (\ref{kkk}) and (\ref{jj}), we have
$$
P( h \le \max\{ U_{[i+1]} - U_{[i]}, i = 1, ..., n-1 \}) \to 0
$$
as $ n \to \infty $. Similarly, we can show that
\beginn
P(\hbar \le  \max\{ V_{[i+1]} - V_{[i]}, i = 1, ..., n-1\}) \to 0
\endn
as $ n \to \infty $. It follows from (\ref{kdf}) and the two equations above that
$$
P(h>
 \max\{ U_{[i+1]} - U_{[i]}, i = 1, ..., n-1 \}, \  \hbar >  \max\{ V_{[i+1]} - V_{[i]}, i = 1, ..., n-1\}) \to 1
$$
as $ n \to \infty$. By Remark 1 and (\ref{cc}), the algorithm converges in probability as $ n \to \infty$.  $\hspace{\fill}\Box$

\

\noindent {\bf Acknowledgements:} The author thanks an associate
editor, a referee and  Professor Z. D. Bai  for their
very valuable comments. The research was partially supported by the National Natural Science Foundation of China (Grant no. 10471061). \bigskip

\section*{References}

\def\ditem{\vspace{-0.2cm} \item}

\begin{description}

%\ditem Ando, T. (1987) Totally positive matrices. \textit{Linear
%Algebra Appl.} \textbf{90}, 165-219.

\ditem Buja, A., Hastie, T. and Tibshirani, R. (1989). Linear
smoothers and additive models (with discussion). \textit{Ann.
Statist.} \textbf{17} 453--555.

\ditem David, H. A. and Nagaraja, H. N. (2003). \textit{Order Statistics}. Wiley, New Jersey.

\ditem Linton, O. and Nielsen, J. P. (1995). A kernel method of
estimating structured nonparametric regression based on marginal
integration. \textit{Biometrika} \textbf{82} 93--100.

\ditem Mammen, E., Linton, O. and Nielsen, J. P. (1999). The
existence and asymptotic properties of a backfitting projection
algorithm under weak conditions. \textit{Ann. Statist.} \textbf{27}
1443--1490.

\ditem Minc, H. (1988). \textit{Nonnegative Matrices}. New York:
Wiley.

\ditem Opsomer, J. D.  and Ruppert, D (1997). Fitting a bivariate
additive model by local polynomial regression. \textit{Ann.
Statist.} \textbf{25} 186--211.

\ditem Romanovsky, V. I. (1970). \textit{Discrete Markov Chains.}
Wolters-Noordhoff Publishing, Groningen, Netherlands.

\ditem Tj{\o}stheim, D. and Auestad, B. (1994). Nonparametric
identification of nonlinear time series: Projections. \textit{J.
Amer. Statist. Assoc.} \textbf{ 89} 1398--1409.

\ditem Wang, L. and Yang, L. (2007). Spline-backfitted kernel
smoothing of nonlinear additive autoregression model. \textit{Ann.
Statist.} \textbf{35} 2474--2503

\end{description}

\end{document}